\documentclass[12pt, oneside]{article}

\usepackage{latexsym}
\usepackage{amsfonts}
\usepackage{amsmath}
\usepackage{amssymb}
\usepackage{amsthm}
\usepackage{xypic}
\usepackage{mathrsfs}
\usepackage{etoolbox}
\usepackage[margin=1in]{geometry}
\usepackage[round]{natbib}

\title{The Witt Ring of a Smooth Projective Curve over a Finite Field}
\author{Funk, Jeanne M. \and Hoobler, Raymond T.}
\date{}

\input xy
\xyoption{all}

\newtheorem{prop}{Proposition}

\newtheorem{thm}{Theorem}

\newtheorem*{pf}{Proof}

\AtBeginEnvironment{pf}{}

\begin{document}

\maketitle

\begin{abstract}

In this paper we calculate the Witt ring $W(C)$ of a smooth geometrically connected projective curve $C$ over a finite field with characteristic other than $2$.  We view $W(C)$ as a subring of $W(k(C))$ where $k(C)$ is the function field of $C$.  The calculation is then completed using classical results for bilinear spaces over fields.
\end{abstract}

Let $k$ be a finite field of characteristic $q\neq 2$. We will show that the Brauer group of a curve $C$ over $k$ vanishes.
The vanishing of the Witt invariant then allows us to represent
any symmetric space in $W(C)$ by a form of rank one or two and allows us to
write out the multiplication and addition table for $W(C)$ and
then recognize it as a quotient of $W(k)[{ }_2Pic(C)]$.

For bilinear spaces over fields, we adopt the notation of~\citep{lam}, in which $<t>$ denotes the rank $1$ space with fixed generator $e$ and whose form takes $(e,e)$ to $t\in k$ and $<t_1,\ldots,t_n>$ is the orthogonal sum $<t_1>\perp\cdots\perp<t_n>$.

The Witt ring $W(k)$ of $k$ is a four element ring.  It consists of $0$, two rank one represented elements $<1>$ (the multiplicative identity) and $<s>$,  and a nontrivial even rank element.  The nontrivial even rank element is $<1,s>$ when $q=1(mod 4)$ and $<1,1>$ when\\
 $q=3(mod 4)$.

$W(k)$ has a few properties which are useful for calculational purposes.

\begin{itemize}
    \item $<1,1>=<s,s>$ for all nondyadic finite fields.
    \item $<1,1,1>=<1>=<-1>$ when $q=1(mod 4)$.
    \item $<1,1,1>=<s>=<-1>$ when $q=3(mod 4)$.
    \item $<1,1,1,1>=0$ for all nondyadic finite fields.
\end{itemize}

Let $C$ be a smooth geometrically connected projective curve over $k$ with generic point $\eta$.  Assume also that $C$ contains a $k$-rational point.

The natural map $W(C)\to W(k(C)):E \to E_{\eta}$ is injective~\citep{balmer} and embeds the Witt ring of $C$ into the Witt ring of its function field.  This suggests the following notation for spaces represented by orthogonal sums of rank $1$ spaces.

 $<\mathscr{L}_{\xi}>$ will denote the Witt class represented by a form which maps to $<\xi>\in W(k(C))$ whose underlying vector bundle is the line bundle $\mathscr{L}$.
$<\mathscr{L}_{1,\xi_1},\ldots,\mathscr{L}_{n,\xi_n}>$ will denote the orthogonal sum $<\mathscr{L}_{1,\xi_1}>\perp\ldots\perp<\mathscr{L}_{n,\xi_n}>$.

Due to the rational point, the ring map $W(k)\hookrightarrow W(C)$
induced by the structure map $C\to Spec(k)$ is also injective.  This
map identifies $<1>$ and $<s>$ with the Witt classes represented by
forms on the structure sheaf.

The following two propositions show that there are $2n$ distinct Witt classes represented by rank $1$ elements where $n$ is the cardinality of the order $2$ Picard group ${ }_2Pic(C)$.

\begin{prop}

Given $\mathscr{L},\mathscr{M}\in{ }_2Pic(C)$, if $<\mathscr{L}_{\xi}>=<\mathscr{M}_{\zeta}>\in W(C)$ then $\mathscr{L}=\mathscr{M}\in{ }_2Pic(C)$ and $<\xi>=<\zeta>\in W(k(C))$.
\end{prop}

\begin{pf}

That $<\xi>=<\zeta>\in W(k(C))$ is trivial.

Let $\mathscr{L}, \mathscr{M}\in {}_2Pic(C)$ and $<\mathscr{L}_{\xi}>=<\mathscr{M}_{\zeta}>\in W(C)$.  Then \\
$<\mathscr{L}_{\xi}>\perp M=<\mathscr{M}_{\zeta}>\perp M'\in Bil(C)$\footnote{$Bil(C)$ is the semiring consisting of isomorphism classes of bilinear spaces on $C$} where $M$ and $M'$ are metabolic spaces of equal rank $2m$.

Taking determinants~\citep{knebusch}  on both sides, we see that\\
 $det(<\mathscr{L}_{\xi}>)det(M)=det(<\mathscr{M}_{\zeta}>)det(M')$, $<\mathscr{L}_{\xi}>(<-1>^m)=<\mathscr{M}_{\zeta}>(<-1>^m)$, and $<\mathscr{L}_{\xi}>=<\mathscr{M}_{\zeta}>\in Bil(C)$. %{<-1>^m is rank1, hence invertable}.
\end{pf}

\begin{prop}

Given $\mathscr{L}\in{ }_2Pic(C)$, there are two Witt classes represented by forms whose underlying space is $\mathscr{L}$.
\end{prop}

\begin{pf}

First, we note that every order $2$ line bundle $\mathscr{L}$ is
equipped with an isomorphism
$\mathscr{L}\otimes\mathscr{L}\to\mathcal{O}_C$, which may be used
to define a nondegenerate symmetric bilinear form on $\mathscr{L}$.

Two rank $1$ forms $(\mathscr{L},\varphi)$ and $(\mathscr{L},\psi)$ may differ at most by a global endomorphism of $\mathscr{L}$ as shown in the following diagram:

$\xymatrix{\mathscr{L}\ar[r]^{\varphi} \ar[d]_{\varphi\psi^{-1}} & \mathscr{L}^{\vee} \ar@{=}[d] \\
            \mathscr{L} \ar[r]_{\psi} & \mathscr{L}^{\vee}}$\\

As $\mathscr{E}nd(\mathscr{L})\cong\mathcal{O}_C$, the global endomorphisms of $\mathscr{L}$ are precisely the units of $k$.  Thus, $\varphi, \psi$ differ by multiplication $m_{\ell}$ by some unit $\ell\in k$.

$\varphi$ and $\psi$ represent the same Witt class precisely when $\varphi = m_{\ell}\circ\psi\circ m_{\ell}^{\vee}$ so that  $\varphi=\ell^2\psi$.

This means that the Witt classes of $C$ associated to $\mathscr{L}$ are in one-to-one correspondence with the square classes of $k$, which correspond to the rank $1$ elements of $W(k)$.
\end{pf}

We note that multiplication by $<s>\neq<1>\in W(k)$ exchanges the two classes of $W(C)$ represented by forms on $\mathscr{L}$.

Since the function field of $C$ is a $C_2$ field we may apply the following classical result for fields ~\citep{lam} Proposition V.3.25.

\begin{prop}
\label{lamthm}
Suppose every form of dimension $5$ over a field $F$ is isotropic.  Then two bilinear spaces $E,E'$ are isometric iff $rk(E)=rk(E')\in\mathbb{Z}$,
$d_{\pm}(E)=d_{\pm}(E')\in W(F)$ and $c(E)=c(E')\in Br(F)$ where $d_{\pm}(-)$ is the signed discriminant and the Witt invariant $c(-)$ is the class of the Clifford algebra in the Brauer group $Br(-)$.
\end{prop}
% (~\citep{lam} pp.119)

The signed discriminant of $E$ is the rank $1$ form $(-1)^{\frac{n(n+1)}{2}}\bigwedge\limits^{n} E$ where $n$ is the rank of $E$.  Note that the signed discriminant of a form $<a_1,\ldots,a_n>$ is $(-1)^{\frac{n(n+1)}{2}}<a_1\cdots a_n>$.

Details regarding Clifford algebras for forms over fields may be
found in~\citep{lam} Chapter V.  A similar definition is used for
bilinear spaces on schemes.  In the current situation it suffices to
know the following:

\begin{itemize}
    \item $c(E)\in Br(C)$.
    \item The natural map $Br(C)\to Br(k(C)):A\mapsto A_{\eta}$ is injective.
    \item $c(E)_{\eta}=c(E_{\eta})$.
\end{itemize}

The following theorem now shows that every bilinear space over $C$ has trivial Witt invariant over $C$ and $k(C)$.

\begin{thm}
\label{Br(C)=0}

$Br(C)=0$
\end{thm}

\begin{pf}

We calculate the cohomological Brauer group $H^2(C_{et}, \mathbb{G}_m)$.

%A finite \'{e}tale map $X'\to X$ between connected varieties said to be a finite Galois covering provided $Gal(X'/X)=|Aut_XX'=deg(X'/X)$.

%Given a finite separable extension $K$ of $k$, $X_K\to X$ is a finite Galois covering with Galois group $G=Gal(K/k)$.

%There is a Hochschild-Serre spectral sequence (~\citep{milne} Theorem III.2.20) $H^p(G,H^q(C_{K,et},\mathbb{G}_m))\Rightarrow H^{p+q}(C_{et},\mathbb{G}_m)$ for each finite separable extension of the base field.

%As cohomology commutes with inverse limits (~\citep{milne} III.1.16),

There is a spectral sequence $H^p(\bar{G},H^q(\bar{C},\mathbb{G}_m))\Rightarrow H^{p+q}(C_{et},\mathbb{G}_m)$~\citep{milne} III.1.16,2.20 where $\bar{C}$ is the extension of $C$ to the separable closure $k_{sep}$ of the base field and \\
$\bar{G}=Gal(k_{sep}/k)$ .  Since $C$ is a curve over a finite field, $\bar{C}$ is also the extension to the algebraic closure.

This is a first quadrant spectral sequence with $E^2$ terms as follows:

$E^2_{0,2}=H^0(\bar{G},H^2(\bar{C},\mathbb{G}_m))$ and $H^2(\bar{C},\mathbb{G}_m)=Br(\bar{C})=0$ due to the fact that $k(\bar{C})$ is a $C_1$ field (Tsen's theorem).  Thus, $E^2_{0,2}=0$.

$E^2_{2,0}=H^2(\bar{G},H^0(\bar{C},\mathbb{G}_m)) =H^2(\bar{G},k_{sep}^{\times})=Br(k)=0$

$E^2_{1,1}=H^1(\bar{G},H^1(\bar{C},\mathbb{G}_m))=H^1(\bar{G},Pic(\bar{C}))$.  There is a short exact sequence
$$\xymatrix{0 \ar[r] & Pic^0(\bar{C}) \ar[r] & Pic(\bar{C}) \ar[r]^{deg} & \mathbb{Z} \ar[r] & 0}$$

Furthermore,  $H^1(\bar{G},Pic^0(\bar{C}))=0$ (Lang's Theorem) and \\
$H^1(\bar{G},\mathbb{Z})=Hom_{cont}(\bar{G},\mathbb{Z})=0$ so that $E^2_{1,1}=H^1(\bar{G},Pic(\bar{C}))=0.$

This shows that the cohomological Brauer group, hence the Brauer group, is trivial.
\end{pf}

Theorem~\ref{Br(C)=0} gives relations $<\mathscr{L}_{\xi},\mathscr{M}_{\zeta}>=<1,\mathscr{L}\mathscr{M}_{\xi\zeta}>$ in $W(C)$ as well as allowing us to show the following.

\begin{prop}

Every class $E\in W(C)$ has a representative of the form $<\mathscr{L}_{\xi}>$ or\\
 $<1,-\mathscr{L}_{\xi}>$ where $<\mathscr{L}_{\xi}>=d_{\pm}(E)$.
\end{prop}

\begin{pf}

Consider the image $E_{\eta}$  of $E$ in the function field.

$E_{\eta}$ has an anisotropic representative $V$ of rank less than five.  We will show that this representative must, in fact, have rank less than three.

We first note that the Witt invariant $c(V)$ is trivial.

The tertiary case now follows from~\citep{lam} Proposition V.3.22.

It remains to show that a rank four form cannot be anisotropic.  We
consider a diagonalisation $<a_1,a_2,a_3,a_4>$ of $V$.  Since $c(V)$ is trivial, Proposition~\ref{lamthm} shows that $V$ is
isometric to $<1,1,1,a_1a_2a_3a_4>=<-1,a_1a_2a_3a_4>\in W(k(C))$.

\end{pf}

Calculations such as $<1,-\mathscr{L}_{\xi}>\perp<\mathscr{M}_{\zeta}>=<1,1,-\mathscr{L}\mathscr{M}_{\xi\zeta}>=<\mathscr{L}\mathscr{M}_{\xi\zeta}>$ now give us the following tables of arithmetic for $W(C)$.

Multiplication:

\begin{tabular}{c|cc}
 & $<\mathscr{L}_{\xi}>$ & $<1,-\mathscr{L}'_{\xi'}>$\\
\hline
$<\mathscr{M}_{\zeta}>$ & $<\mathscr{L}\mathscr{M}_{\xi\zeta}>$ & $<1,-\mathscr{L}'_{\xi'}> $\\
$<1,-\mathscr{M}'_{\zeta'}>$ & $<1,-\mathscr{M}'_{\zeta'}>$& $0$
\end{tabular}\\

Note, in particular, that $<\mathscr{L}_{\xi}><\mathscr{L}_{\xi}>=<1>$ and\\
 $<\mathscr{L}_{\xi}><\mathscr{L}_{s\xi}>=<s>$

Addition:

\begin{tabular}{c|cc}
 & $<\mathscr{L}_{\xi}>$ & $<1,-\mathscr{L}'_{\xi'}>$\\
\hline
$<\mathscr{M}_{\zeta}>$ & $<1,\mathscr{L}\mathscr{M}_{\xi\zeta}>$ & $<\mathscr{M}\mathscr{L}'_{\zeta\xi'}> $\\
$<1,-\mathscr{M}'_{\zeta'}>$ & $<\mathscr{M}'\mathscr{L}_{\zeta'\xi}>$& $<1,-\mathscr{M}'\mathscr{L}_{\zeta'\xi}>$
\end{tabular}\\

Note that $<\mathscr{L}_{\xi}>\perp<\mathscr{L}_{\xi}>=<1,1>$ and $<\mathscr{L}_{\xi}>\perp<\mathscr{L}_{s\xi}>=<1,s>$.

In fact, $W(C)$ can be expressed quite nicely as a quotient of the group ring $W(k)[{ }_2Pic(C)]$.

\begin{thm}

$W(C)\cong W(k)[{ }_2Pic(C)]/\mathscr{R}$ where $\mathscr{R}$ is generated by the relations of the form
 $$<1>-<u>\mathscr{L}-<v>\mathscr{M}+<uv>\mathscr{L}\mathscr{M}$$
 with $<u>,<v>\in W(k)$ and $\mathscr{L},\mathscr{M}\in{ }_2Pic(C)$.
\end{thm}

\begin{pf}

Fix a form $<\mathscr{L}_{\xi}>$ on each $\mathscr{L}\in{ }_2Pic(C)$.

Define a map $W(k)[{}_2Pic(C)]/\mathscr{R}\to W(C)$ by sending $<s>\mathscr{L}$ to $<\mathscr{L}_{s\xi}>$ and extending by linearity.

This map is clearly a well defined surjection of commutative rings.

To show injectivity let $f\in W(k)[{}_2Pic(C)]/\mathscr{R}$ map to a form \\
$E=<\mathscr{L}_{1,u_1\xi_1},\ldots,\mathscr{L}_{n,u_n\xi_n}>=0 \in W(C)$.  Using the relations $<1,1>=<-1,-1>$ and\\
 $<\mathscr{L}_{u\xi},\mathscr{M}_{v\zeta}>=<1,\mathscr{L}\mathscr{M}_{uv\xi\zeta}>$ we may rewrite $E$ as $<\pm 1,\mathscr{L}_1\cdots\mathscr{L}_{n,u_1\cdots u_n\xi_1\cdots\xi_n}>$.   Using the corresponding relations $<1>+<1>+<1>+<1>$ and\\
 $<1>-<u>\mathscr{L}-<v>\mathscr{M}+<uv>\mathscr{L}\mathscr{M}$ in $W(k)[{}_2Pic(C)]/\mathscr{R}$, we may write $f$ as $<\pm 1>+<u_1\cdots u_n>\mathscr{L}_1\cdots\mathscr{L}_n$.  In the $+1$ case, $E$ and $f$ are trivial precisely when $<u_1\cdots u_n\mathscr{L}_1\cdots\mathscr{L}_n>=<1>$.  In the $-1$ case, both $E$ and $f$ are trivial precisely when $<u_1\cdots u_n\mathscr{L}_1\cdots\mathscr{L}_n>=<-1>$.
\end{pf}

\bibliographystyle{plainnat}
\bibliography{bib1}

\end{document}